\documentclass[11pt,a4paper,leqno]{article}
\usepackage{amsmath}
\usepackage{amssymb}
\advance\topmargin by -2.2cm
\newcommand{\bbr}{I\!\!R}

\newcommand{\bbn}{I\!\!N}
\newcommand{\bbz}{Z\!\!\!Z}
\newcommand{\bbf}{I\!\!F}
\newcommand{\bbg}{I\!\!\!\!G}
 
\newcommand{\cala}{{\cal A}}
\newcommand{\calb}{{\cal B}}
\newcommand{\calc}{{\cal C}}

\newcommand{\cale}{{\cal E}}

\newcommand{\calg}{{\cal G}}

\newcommand{\call}{{\cal L}}

\newcommand{\caln}{{\cal N}}

\newcommand{\barr}{\begin{array}}
\newcommand{\earr}{\end{array}}
\newcommand{\beqq}{\begin{equation}}
\newcommand{\eeqq}{\end{equation}}
\newcommand{\beao}{\begin{eqnarray*}}
\newcommand{\eeao}{\end{eqnarray*}\noindent}
\newcommand{\beam}{\begin{eqnarray}}
\newcommand{\eeam}{\end{eqnarray}\noindent}

\newcommand{\halmos}{\quad\hfill\mbox{$\Box$}}

\newcommand{\la}{\lambda}

\newcommand{\si}{\sigma}
\newcommand{\al}{\alpha}
\newcommand{\vth}{\vartheta}

\newcommand{\Om}{\Omega}

\newcommand{\vep}{\varepsilon}

\newcommand{\ov}{\overline}
\newcommand{\wh}{\widehat}
\newcommand{\wt}{\widetilde}

\newcommand{\ul}{\underline}

\newcommand{\lra}{\longrightarrow}

\newcommand{\nto}{n\to\infty}

\newcommand{\lala}{\lambda\!\!\lambda}

\begin{document}

{\huge 
On LAN for parametrized continuous periodic signals\\ in a time inhomogeneous diffusion}
\vskip0.5cm
{\bf R.\ H\"opfner}, Johannes Gutenberg Universit\"at Mainz\\
{\bf Yu.\ Kutoyants}, Universit\'e du Maine, Le Mans\\

{\bf Abstract: } We consider a diffusion $(\xi_t)_{t\ge 0}$ whose drift involves some $T$-periodic signal. $T$ is fixed and known, whereas the signal depends on a $d$-dimensional parameter $\vth\in\Theta$.  Assuming positive Harris recurrence of the grid chain $(\xi_{kT})_{k\in\bbn_0}$ and exploiting the periodic structure in the semigroup, we work with path segments and limit theorems for certain functionals (more general than additive functionals) of the process to prove local asymptotic normality (LAN). Then we consider several estimators for the unknown parameter.  
{\bf MSC: \quad 62 F 12 , 60 J 60} \\  

\today\\


We consider a problem of parameter estimation in a time-inhomogenous diffusion $(\xi_t)_{t\ge 0}$ whose semigroup has a $T$-periodic structure, and whose drift involves a signal --deterministic, continuous, of known periodicity $T$-- which is  parametrized by some $d$-dimensional parameter $\vth\in\Theta$. Our main assumption on the process is positive Harris recurrence of the grid chain $(\xi_{kT})_{k\in\bbn_0}$ which implies (as in [HK 10]) positive Harris recurrence of the path-segment chain 
$$
\mathbb{X} = \left( \mathbb{X}_k\right)_k  \quad,\quad    \mathbb{X}_k \;:=\; \left( \xi_{(k-1)T+s}Ê\right)_{0\le s\le T}  
$$
taking values in $C([0,T])$: this allows for limit theorems using the path-segment chain. Our main statistical assumption is that the derivative of the signal with respect to the parameter is again a $T$-periodic function, and that the parametrization is sufficiently smooth, in a suitable $L^2$-sense related to the periodic structure of the semigroup.   

We prove local asymptotic normality (LAN) for the sequence of statistical models corresponding to observation of a trajectory of $\xi$ --~with unknown $\vth$~-- over a long time interval $[0,nT]$, and state a local asymptotic minimax theorem. See LeCam ([LC 68], [LC 90]), H\'ajek [H 71], Ibragimov and Khasminskii [IH 81] for background on LAN and convergence of estimators, see also the survey paper on local asymptotics by Davies [D 85], and see Liptser and Shiryaev [LS 81], Jacod and Shiryaev [JS~87], Kutoyants [K 04] on background on statistical models and likelihood ratio processes in diffusions. 'Signal in white noise' problems seem to appear first in Ibragimov and Khasminskii [IH 81]. The limit theorems which we use for time-inhomogeneous diffusions $(\xi_t)_{t\ge 0}$ with periodic structure  in the semigroup are from [HK 10, Section 2]. 

We then consider a sequence of minimum distance estimators (MDE) which via LeCam's one step correction can be improved to get an asymptotically efficient estimator sequence, thus providing an easy-to-calculate and explicit alternative to  
maximum likelihood estimators and Bayes estimators which under some conditions are asymptotically efficient for the unknown parameter.  \\

\section*{1.~ Local asymptotic normality}

Our setting is as follows. The observed diffusion is inhomogeneous in time 
\beqq\label{process}
d\xi_t  \;=\;  \left[ S(\vth,t) + b(\xi_t) \right]dt   \;+\;  \si(\xi_t)\,dW_t  \quad,\quad t\ge 0 \;,\; \xi_0=x_0  
\eeqq
and depends on an unknown parameter $\vth$ ranging over an open set $\Theta\subset\bbr^d$. The drift in (\ref{process}) involves a deterministic signal $S(\cdot,\cdot)$, which is a continuous function $\Theta{\times}[0,\infty) \to \bbr$ satisfying 
\beqq\label{Vper}
\mbox{for every $\vth\in\Theta$ fixed: }\quad
S(\vth,\cdot)  :   [0,\infty)\to\bbr  \quad\mbox{is  $T$-periodic}   
\eeqq
for some period $T$ (fixed and known, and not depending on $\vth$). We write throughout the paper 
$$
 i_T(t) \;:=\; t \;\mbox{modulo}\; T   \;. 
$$
The functions $b(\cdot)$ and  $\si(\cdot)$ are assumed Lipschitz, and $\si(\cdot)$ strictly positive. Hence for all values of $\vth\in\Theta$,  we have Lipschitz and linear growth conditions for the time-dependent coefficients of equation (\ref{process}), and thus existence  and pathwise uniqueness for its solution ([S 65], [KS 91]). We discuss 'ergodicity properties' for the time inhomogeneous diffusion $\xi$. Write $(P^{(\vth)}_{s,t})_{0\le s<t<\infty}$ for the semigroup of  the process (\ref{process}). As a consequence of (\ref{Vper}), this semigroup is $T$-periodic in the sense 
\beqq\label{periodicsemigroup}
P^{(\vth)}_{s,t}(x,dy) \;=\; P^{(\vth)}_{s+kT,t+kT}(x,dy) \quad\mbox{for all $k\in\bbn_0$ and all $0\le s<t<\infty$} \;. 
\eeqq
We  assume 
$$
\begin{array}{l}
\mbox{ for every $\vth\in\Theta$, the grid chain $(\xi_{kT})_{k\in\bbn_0}$ is positive recurrent }\\
\mbox{ in the sense of Harris with invariant probability $\mu^{(\vth)}$ on $(\bbr,\calb(\bbr))$ } \\[2mm]
\end{array}
\leqno{(H1)}
$$
and decompose the path of $\xi$ into a Markov chain of $T$--segments 
$$
\mathbb{X} = (\mathbb{X}_k)_k  \quad\mbox{defined by}\quad  \mathbb{X}_k \;:=\; (\xi_{(k-1)T+s})_{0\le s\le T} \;,\; k\ge 1 
$$
which takes values in the path space $(C_T,\calc_T)$ of continuous functions $[0,T]\to\bbr$, using as initial value $\mathbb{X}_0$ some $\al\in C_T$ such that $\al(T)=x_0$. As a consequence of  (\ref{periodicsemigroup}), the chain $\mathbb{X}$ is time homogeneous. As a consequence of  $(H1)$, $\,\mathbb{X}$ will be positive recurrent in the sense of Harris under $\vth$ (see [HK 09, theorem 2.1]), with invariant probability  $m^{(\vth)}$ on $(C_T,\calc_T)$ whose finite-dimensional distributions are given in terms of $\mu^{(\vth)}$ and $P^{(\vth)}_{s,t}$, $0\le s<t\le T$.  
In order to obtain limit theorems for log-likelihoods related to observation of the time-inhomogeneous process $\xi$ over a long time interval, for asymptotics of local models at $\vth$, or for rescaled estimation errors of interesting estimators, we shall always exploit the ergodicity of the time homogeneous $C_T$-valued chain $\mathbb{X}$. \\

We turn to properties of the parametrization and list the properties of the parametrization which we shall need. Our identifiability  condition is 
$$
\mbox{ for all $\vth\in\Theta$,  $\vep>0$:  }\quad 
\inf\limits_{ |\zeta - \vth| > \vep}\; \sup\limits_{0\le s\le T}\; \left| S(\zeta,s)-S(\vth,s)\right| \;>\; 0 \;. 
\leqno{(H2)}
$$
A sufficiently general  differentiability condition is as follows.  Restrict the continuous periodic functions $\{ S(\zeta,\cdot) : \zeta\in\Theta \}$  to $[0,T]$,  assume that all measures 
\beqq\label{measureslambda}
 \la^{(\vth)}(ds) \;:=\;  \left[ \mu^{(\vth)}\! P^{(\vth)}_{0,s} (\frac{1}{\si^2}) \right] ds \quad\mbox{on}\;\; ([0,T],\calb([0,T])) 
 \quad,\quad \vth\in\Theta  
\eeqq
are finite,  and consider spaces 
$$
\mathbb{H}^{(\vth)} \;=\; L^2( [0,T], \calb([0,T]),  \la^{(\vth)} )  \quad,\quad \vth\in\Theta \;.  
$$
Assume that for every $\vth\in\Theta$, there is a function  
$$
\dot S(\vth,\cdot) = \left(\begin{array}{l} \dot S_1(\vth,\cdot) \\ \cdots \\ \dot S_d(\vth,\cdot) \end{array}\right) \quad,\quad 
\dot S_j(\vth,\cdot)  \;\in\;  \mathbb{H}^{(\vth)}   \;\;\mbox{for $j=1,\ldots,d$}
$$
such that the following holds: 
$$
\left\{\begin{array}{l}
\rho_\vep(\vth,\cdot) := \sup\limits_{|\zeta-\vth|<\vep} \left| \frac{ S(\zeta,\cdot) - S(\vth,\cdot) - (\zeta-\vth)^\top \dot S(\vth,\cdot) }{ \zeta-\vth} \right| 
\quad\in\;\; \mathbb{H}^{(\vth)} \quad\mbox{for some $\vep=\vep(\vth)>0$} \\
\rho_\vep(\vth,\cdot) \;\lra\; 0\quad \mbox{ in $\mathbb{H}^{(\vth) }\;$ as $\;\vep\downarrow 0$}\;. 
\end{array}\right. 
\leqno{(H3)}
$$
This is slightly more than Fr\'ech\'et differentiability of $\zeta\to S(\zeta,\cdot)$ at  $\zeta=\vth$ in the Hilbert space $\mathbb{H}^{(\vth)}$. The assumption that the derivative with respect to the parameter $t\to\dot S(\vth,t)$  is again a $T$-periodic function --the key to the results presented here-- is a strong assumption;  e.g., it rules out  estimation of an unknown periodicity (`frequency modulation`) as considered in [IH 81, p.\ 209] or [CLM 06]. We call the $d{\times}d$ matrix  
\beqq\label{fisherinformation}
I^{(\vth)} \;:=\; \int_0^T \dot S(\vth,s)\,  \dot S(\vth,s)^{\!\top} \la^{(\vth)}(ds)   
\eeqq
Fisher information at $\vth$, and add the assumption 
$$
\mbox{$I^{(\vth)}$ is invertible for all $\vth\in\Theta$}\;.  
\leqno{(H4)}
$$  
We introduce two additional assumptions. For $S := \{ \al 2^{-k} : k\in\bbn_0\,,\, \al\in\bbz^d \}$ the set of dyadic numbers in $\bbr^d$, we shall need 
$$
\left\{\begin{array}{l}
\wt\rho_\vep(\vth,\cdot) := \sup\limits_{\zeta \in S\cap\Theta , |\zeta-\vth|<\vep} \left|  \dot S(\zeta,\cdot) - \dot S(\vth,\cdot) \right| 
\quad\in\;\; \mathbb{H}^{(\vth)} \quad\mbox{for some $\vep=\vep(\vth)>0$} \\
\wt\rho_\vep(\vth,\cdot) \;\lra\; 0\quad \mbox{ in $\mathbb{H}^{(\vth) }\;$ as $\;\vep\downarrow 0$} 
\end{array}\right. 
\leqno{(H5)}
$$
for all $\vth\in\Theta$ together with 
$$
\sup_{\zeta\in S\cap\Theta \,,\, |\zeta-\vth|<\vep}\; \left| I^{(\zeta)} - I^{(\vth)} \right| \;\;\lra\;\;   0  
\quad\mbox{as}\quad \vep\,\downarrow\,0  \;. 
\leqno(H6)
$$

\vskip0.8cm
Let $Q^\vth$ denote the law of the solution of (\ref{process})+(\ref{Vper}) under $\vth\in\Theta$, a law on the canonical path space $(C,\calc)$ of continuous functions $[0,\infty)\to\bbr$ equipped with the metric of locally uniform convergence. Write $Q^\vth_t$ for the restriction of $Q^\vth$ to the $\si$-field $\calg_t$ of events before time $t^+$, and $\bbg = (\calg_t)_{t\ge 0}$. The likelihood process $L^{\zeta / \vth}$ of $Q^\zeta$ with respect to $Q^\vth$ relative to $\bbg$  (cf. [LS 81], [JS 87], [K 04]) is given by 
\beqq\label{likelihood}
L^{\zeta / \vth} 
\;=\;  \cale_\vth\left( \int_0^\cdot \frac{S(\zeta,t)-S(\theta,t)}{\si^2(\eta_t)}\; dm^{(\vth)}_t \right) 
\;=\;  \cale_\vth\left( \int_0^\cdot \frac{S(\zeta,t)-S(\theta,t)}{\si(\eta_t)}\; dB^{(\vth)}_t \right) 
\eeqq
where $\eta=(\eta_t)_{t\ge 0}$ denotes the canonical process on $(C,\calc, \bbg)$,  $\;m^{(\vth)}$ its $(\bbg,Q^\vth)$-martingale part, and where $B^{(\vth)} := \int_0^\cdot \frac{1}{\si(\eta_s)}dm^{(\vth)}_s$ is a  $(Q^\vth, \bbg)$-Brownian motion. The following is a '2nd Le Cam lemma' for time inhomogeneous diffusions of type (\ref{process})+(\ref{Vper}):\\

{\bf 1.1 Theorem : } Under ($H1$)+($H3$), in the sequence of experiments 
$$
\left(\, C \,,\, \calg_{nT} \,,\, \left\{ Q^{\zeta}_{nT} : \zeta\in\Theta \right\} \,\right) \quad,\quad \nto \;, 
$$
we have LAN at every point $\vth\in \Theta$,  with local scale $n^{-1/2}$ and Fisher information 
$I^{(\vth)}$ given by (\ref{fisherinformation}). More precisely, for every $\vth\in\Theta$ and arbitrary bounded sequences $(h_n)_n$ in $\bbr^d$, we have a quadratic decomposition of log-likelihood ratios  
\beqq\label{quadraticdecomposition}
\log L_{nT}^{(\vth + n^{-1/2} h_n) / \vth}  \;=\; h_n^\top \Delta_n^{(\vth)} \;-\; \frac12 h_n^\top I^{(\vth)} h_n \;+\; o_{Q^\vth}(1) \quad\mbox{as $\nto$} 
\eeqq
with score 
$$
\Delta_n^{(\vth)} \;=\; \frac{1}{\sqrt{n}} \int_0^{nT} \frac{\dot S(\vth,s)}{\si(\eta_s)}\, dB^{(\vth)}_s 
$$
such that  
\beqq\label{convergence of score}
\call (\, \Delta_n^{(\vth)} \mid Q^\vth  )  \;\lra\;  \caln (\, 0 \,,\, I^{(\vth)}  \,) \quad\mbox{(weak convergence in $\bbr^d$, $\nto$)} \;. 
\eeqq
 
\vskip0.8cm
The proof of theorem 1.1 will be given in section 2 below (see 2.1). For background on LAN see [L 69], [H 70], [IH 81], [D 85], [LY 90], [K 04], or also [H 08] for some more details. A well-known and powerful consequence of LAN is the local asymptotic minimax theorem (corollary 1.1' below, from [L~90, p.\ 83]; a result using neighbourhoods in a different way is [IH 81], with additional uniformity assumptions which we do not make here). Fix any reference point $\vth\in\Theta$, write 
\beqq\label{defcentralsequence}
Z_n^{(\vth)} \;:=\; (I^{(\vth)})^{-1}\Delta_n^{(\vth)} \;\;,\;\; n\ge 1
\eeqq
for the central sequence in the sequence of local models 
$$
\cale_n(\vth)  \;:=\;  
\left\{ P_{n,\vth,h} := Q^{\vth + n^{-1/2} h}_{nT} \,:\, h\in\bbr^d  \;\mbox{such that}\;\,  \vth{+}n^{-1/2}h \,\in\,\Theta \right\} \;\;,\;\; n\ge 1
$$
at $\vth$, and let $\ell(\cdot)$ denote any loss function which is continuous, subconvex and bounded. Then \\

{\bf 1.1' Corollary :} (LeCam [L 90, p.\ 83]) Under $(H1)$ and $(H3){+}(H4)$, the following holds.  
For any sequence of $\calg_{nT}$-measurable estimators $\wt\vth_n$, for any $\vth\in\Theta$,  
$$
\sup_{c<\infty}\; \liminf_{\nto}\; \sup_{|h|\le c}\; 
E_{(\vth + n^{-1/2} h)}\left( \ell\left(  n^{1/2}( \wt\vth_n - (\vth + n^{-1/2} h) ) \right)  \right) 
\quad\ge\quad 
\int_{-\infty}^\infty \ell(z)\, \caln(0,(I^{(\vth)})^{-1})(dz) \;. 
$$
Sequences of $\calg_{nT}$-measurable estimators $\vth^*_n$ with the property 
\beqq\label{criterionattheta}
n^{1/2}( \vth^*_n - \vth )  \;=\;  Z_n^{(\vth)} 
\;+\; o_{Q^\vth}(1) \quad\mbox{as $\nto$} 
\eeqq
attain the local asymptotic minimax bound at $\vth$, and have 
$$
\liminf_{\nto}\; \sup_{|h|\le c}\; 
E_{(\vth + n^{-1/2} h)}\left( \ell\left(  n^{1/2}( \vth^*_n - (\vth + n^{-1/2} h) ) \right)  \right) 
\quad=\quad 
\int_{-\infty}^\infty \ell(z)\, \caln(0,(I^{(\vth)})^{-1})(dz)   
$$
for every $0<c<\infty$. \\

Estimator sequences attaining the local asymptotic minimax bound at $\vth$ are termed efficient at $\vth$. The criterion (\ref{criterionattheta}) characterizes efficiency at $\vth$ in terms of stochastic equivalence under $\vth$ of rescaled estimation errors at $\vth$ and central sequence at $\vth$. Using Le Cam's 'one step correction', efficient estimator sequences can be constructed explicitely provided we dispose of at least one preliminary estimator sequence which is $n^{1/2}$-consistent. The construction which we give here follows Davies [D 85]. \\

{\bf Theorem 1.2 : } Assume $(H1)$, $(H3){+}(H4)$, $(H5){+}(H6)$, and let  $\wt\vth_n$ denote some sequence of $\calg_{nT}$-measurable estimators with the property 
$$
\mbox{for every $\vth\in\Theta$}\;:\quad \call\left( n^{1/2}( \wt\vth_n - \vth ) \mid Q^{(\vth)}\right) \;,\; n\ge 1\;,\; \mbox{is tight in $\bbr^d$}\;.   
$$
For $\al=(\al_1,\ldots,\al_d)\in\bbz^d$ and $k\in\bbn_0$, consider cubes $C(k,\al) := \mathop{\sf X}\limits_{j=1}^d [\al_j 2^{-k} , (\al_j{+}1)2^{-k} [$ in $\bbr^d$. Fixing some point $\vth_0\in\Theta$, we discretize the estimator $\wt\vth_n$ to 
$$
\wt{\wt \vth}_n \;:=\; \vth_0\; +\!\! \sum_{\al\in\bbz^d \;\mbox{s.t.}\; C(n,\al)\subset\Theta} (\al 2^{-n}-\vth_0)\; 1_{C(n,\al)}(\wt \vth_n) 
$$
which takes a countable number of values in $\Theta$. Then  
$$
\vth_n^*  \;:=\;  \wt{\wt \vth}_n + n^{-1/2}\, (I^{(\wt{\wt \vth}_n)})^{-1}\, \Delta_n^{(\wt{\wt \vth}_n)} 
$$
is a sequence of $\calg_{nT}$-measurable estimators for $\vth\in\Theta$ which has the property (\ref{criterionattheta}) at all points $\vth\in\Theta$. \\

Thus by corollary 1.1', the sequence $\left(\vth_n^*\right)_n$ constructed in theorem 1.2 is efficient at all points $\vth\in\Theta$.  Since sequences meeting (\ref{criterionattheta}) are regular at $\vth$ in the sense of H\'ajek ([H 70]), they are also regular and efficient in the sense of H\'ajek's convolution theorem.  \\

{\bf 1.3 Example: } We look to equations (\ref{process}) of Ornstein-Uhlenbeck type
$$
b(x)=-\beta x \quad\mbox{for some}\;\; \beta>0 \quad,\quad \si(\cdot)\equiv\si>0 \quad\mbox{constant} \;. 
$$
1)~ It follows from [HK 10, example 2.3] that $(H1)$ holds for arbitrary continuous and $T$-periodic functions $t\to S(\vth,t)$, $\vth\in\Theta$. The diffusion coefficient being constant, the measure $\la^{(\vth)}$ in 
(\ref{measureslambda}) is Lebesgue measure on $[0,T]$ multiplied by the constant $\frac{1}{\si^2}$. It does not depend on $\vth\in\Theta$. 

2)~ Let $\Theta\subset\bbr^d$ be open, let $(K_m)_m$ denote an increasing sequence of compacts with $K_m\subset int(K_{m+1})$, $m\ge 1$,  and $\Theta = \bigcup_m K_m$. For every $t\in[0,T]$ fixed, let the mapping $\zeta\to S(\zeta,t)$ be $\calc^1$ on $\Theta$. Then  
\beqq\label{sufficientcondition2} 
\mbox{for every $m$ :}\quad 
\sup_{t\in[0,T]}\; \sup_{\zeta\in K_m}\; \sup_{|u|=1}\; \left| u^\top \dot S(\zeta,t) \right| \quad<\quad \infty 
\eeqq
is a sufficient condition implying $(H3)$+$(H5)$+$(H6)$. Under (\ref{sufficientcondition2}), at every point $\vth\in\Theta$,  
$$ 
\sup_{\zeta\in K_m}\; \frac{1}{|\zeta-\vth|^2}\; \int_0^T (S(\zeta,t)-S(\vth,t))^2\, dt \quad<\quad \infty  
$$ 
provided $m$ is large enough for $\vth\in\,${\em int}$(K_m)$, and we may consider the following partially converse condition:  
for all $\vth$, assume that there is some $\vep=\vep(\vth)>0$ such that
\beqq\label{sufficientcondition3}
\inf_{\zeta\in\Theta}\; \frac{1}{|\zeta-\vth|^2\wedge\vep}\; \int_0^T (S(\zeta,t)-S(\vth,t))^2\, dt \quad>\quad 0 \;. 
\eeqq
Then under (\ref{sufficientcondition2}), condition (\ref{sufficientcondition3}) implies both $(H2)$+$(H4)$ together. 
Hence, if  $\zeta\to S(\zeta,t)$ is $\calc^1$ on $\Theta$ for fixed $t$, (\ref{sufficientcondition2})+(\ref{sufficientcondition3}) are sufficient conditions for $(H2)$-$(H6)$.

3)~ Fix some continuous $T$-periodic function $f(\cdot)$, non-constant. Put $\Theta:=(0,T){\times}(0,\infty)$, write $\vth=(\vth_1,\vth_2)$, and define (with notation $i_T(t)$ for $t$ modulo $T$ as above)
$$
S(\vth,t) \;:=\; \vth_2\, f( i_T(t) - \vth_1 ) \quad,\quad  t\in\bbr \;.  
$$
Then (\ref{sufficientcondition2})+(\ref{sufficientcondition3}) and thus  $(H2)$-$(H6)$ hold. By $T$-periodicity of $f(\cdot)$, the Fisher information in (\ref{fisherinformation}) is 
$$
I^{(\vth)} \;=\; \frac{1}{\si^2} \left( \begin{array}{ll} 
\vth_2^2 \int_0^T[f']^2(s) ds & -\vth_2 \int_0^T [f' f](s) ds \\
-\vth_2 \int_0^T [f' f](s) ds & \int_0^T f^2(s) ds 
\end{array}\right)
$$
which does not depend on $\vth_1\in(0,T)$. 

4)~ Assume that $T$ is large compared to $2$, put $\Theta=(0,T)$, and define 
$$
S(\vth,t) \;:=\; f(t-\vth) \quad\mbox{with}\quad  f(x) = 1_{\{|i_T(t)|\le 1\}}(1-|i_T(t)|) \;,\; t\in\bbr \;. 
$$
Then assumptions $(H2)-(H6)$ are satisfied with  
$$
\dot S(\vth,t) \;=\; \left( 1_{(\vth,\vth+1)} - 1_{(\vth-1,\vth)} \right) (i_T(t)) \quad,\quad   t\in\bbr   
$$
and Fisher information 
$$
I^{(\vth)} \;=\; \frac{2}{\si^2} \quad,\quad \vth\in\Theta   
$$
not depending on $\vth\in\Theta$. Similiarly, if we replace the definition of $f$ above by 
$$
f(x) = 1_{\{|i_T(t)|\le 1\}}(1-|i_T(t)|)^\al \;,\; t\in\bbr 
$$
with some $\al>1$, then assumptions $(H2)-(H6)$ are satisfied: we have  
$$
\dot S(\vth,t) \;=\; g(t-\vth) \quad\mbox{with}\quad  
g(t) \;:=\; \al\, (1-|i_T(t)|)^{\al-1}\; \left( 1_{(0,1)} - 1_{(-1,0)} \right)(i_T(t)) 
$$
and the Fisher information does not depend on $\vth\in\Theta$. This is analogous to some examples in [L~90, p.\ 32]: 
$\dot S(\vth,\cdot)$ as element of $\mathbb {H}^{(\vth)}$ is defined for $\al>\frac12$, $(H3)$ holds for $\al>\frac12$ by dominated convergence, but $(H5)$ is violated in case $1>\al>\frac12$. \halmos\\


\vskip1.0cm
\newpage
\section*{2.~ Proofs}

In this section, we assume that the strong Markov process $\xi$ lives on some  $( \Om, \cala, \bbf, (P_x)_{x\in\bbr})$; 'almost surely' means almost surely with respect to every $P_x$, $x\in\bbr$; $\bbf^\xi$ is the filtration generated by $\xi$, and an $\bbf^\xi$--increasing process is an $\bbf^\xi$--adapted c\`adl\`ag process $A = (A_t)_{t\ge 0}$ with nondecreasing paths and $A_0=0$, almost surely. 

$(C_T,\calc_T)$ --~the space of all continuous functions $\al:[0,T]\to\bbr$ with  Borel $\si$-field $\calc_T$~--   is a Polish space, 
and $\calc_T$ is also generated by the coordinate projections $\pi_t$, $0\le t\le T$.  The one-step-transition kernel  $Q(\cdot,\cdot)$ of the chain $\mathbb{X}$ of  $T$-segments in the path of $\xi$ is given  by 
$$
Q(\al,F) \;:=\; P\left(\, (\xi_s)_{0\le s\le T} \in F \mid \xi_0=\al(T)\, \right) \;,\quad \al\in C_T \;,\; F\in\calc_T \;. 
$$
For Harris processes in discrete time, we refer to Revuz [R 75] or Meyn and Tweedie [MT 93]. For Harris processes in continuous time see Az\'ema, Duflo and Revuz [ADR 69] or R\'evuz and Yor [RY~91, Ch.\ 10.3]. The following (with sub- or superscripts $\vth$ suppressed from notation) is a precise statement of the 'ergodicity properties'  mentioned in section 1. \\

{\bf Theorem A: } ( [HK 10, theorem 2.1]) Assume that  the grid chain $(\xi_{kT})_{k\in\bbn_0}$ is positive recurrent in the sense of Harris, and write $\mu$ for its invariant probability on $(\bbr,\calb(\bbr))$.

a) Then the chain $\mathbb{X} = (\mathbb{X}_k)_{k\in\bbn_0}$  of $T$-segments in the path of $\xi$ is positive recurrent in the sense of Harris. Its  invariant probability  is the unique law $m$ on $(C_T,\calc_T)$ such that 
\beqq\label{invmeasureTper-1}
\left\{ \begin{array}{l}
\mbox{for arbitrary $0=t_0<t_1<\ldots<t_{\ell}<t_{\ell+1}=T$ and $A_i\in\calb(\bbr)$} \;, \\
m\left( \left\{ \pi_{t_i}\in A_i \,,\,0\le i\le \ell{+}1 \right\}\right) \quad\mbox{is given by}\\
\int\ldots\int \mu(dx_0)\, 1_{A_0}(x_0)\, \prod_{i=0}^\ell P_{t_i,t_{i+1}}(x_i,dx_{i+1})\, 1_{A_{i+1}}(x_{i+1}) \;. 
\end{array} \right.
\eeqq

b) For every $\bbf^\xi$--increasing process $A = (A_t)_{t\ge 0}$ with the property
\beqq\label{property}
\left\{ \begin{array}{l}
\mbox{there is some function $F:C_T\to\bbr$, nonnegative, $\calc_T$-measurable, }\\
\mbox{satisfying}\quad   m(F) := \int_{C_T} F(\al)\, m(d\al) \;<\;  \infty \;, \;\;\mbox{such that} \\
A_{kT} \;=\; \sum\limits_{j=1}^k F(\mathbb{X}_k)  
\;=\; \sum\limits_{j=1}^k  F\left(\, (\xi_{(k-1)T+s})_{0\le s\le T} \,\right) \;,\; k\ge 1  
\end{array} \right.
\eeqq
we have the strong law of large numbers 
$$
\lim_{t\to\infty}\, \frac{1}{t}\, A_t \;\;=\;\; \frac{1}{T}\, m(F) \quad\mbox{almost surely} \;. 
$$

\vskip0.8cm
This strong law of large numbers will be the key tool to prove theorem 1.1. We use the notations which have been introduced in section 1. \\

{\bf 2.1 Proof of theorem 1.1 : } We assume $(H1)$ and $(H3)$.  

1)~ Fix $\vth$. For $h\in\bbr^d$, consider 
$$
\log L_{nT}^{(\vth+n^{-1/2}h) / \vth } \;=\; \int_0^{nT} \frac{S(\vth+n^{-1/2}h,s)-S(\vth,s)}{\si(\eta_s)}\, dB^{(\vth)}_s 
\;-\; \frac12 \int_0^{nT} \left( \frac{S(\vth+n^{-1/2}h,s)-S(\vth,s)}{\si(\eta_s)} \right)^2 ds 
$$
under $Q^\vth$. Adding $\pm n^{-1/2} h^\top \dot S(\vth,\cdot)$ in the numerators of the integrands,  we separate leading terms
\beqq\label{leadingterm0}
h^\top \Delta^n_T(\vth)  \;-\; \frac12 h^\top I^n_T(\vth)\, h 
\eeqq
defined by  
\beqq\label{leadingterm1}
\Delta^n_t(\vth) \;:=\;  \frac{1}{\sqrt{n}} \int_0^{tn} \frac{\dot S(\vth,s)}{\si(\eta_s)}\, dB^{(\vth)}_s   \;\;,\;\;   t\ge 0 
\eeqq
\beqq\label{leadingterm2}
I^n_t(\vth) \;:=\; \frac1n \int_0^{tn} \left( \frac{ \dot S(\vth,s)\;\dot S(\vth,s)^{\!\top} }{\si^2(\eta_s)} \right) ds  \;\;,\;\;   t\ge 0
\eeqq
from  remainder terms 
\beqq\label{remainderterm0}
N^n_T(h) \;-\; \frac12  U^n_T(h) \;-\; V^n_T(h)
\eeqq
defined by 
\beqq\label{remainderterm1}
N^n_t(h) \;:=\;  \int_0^{tn} \frac{S(\vth+n^{-1/2}h,s)-S(\vth,s) -  n^{-1/2} h^\top \dot S(\vth,s)}{\si(\eta_s)}\, dB^{(\vth)}_s  \;\;,\;\;   t\ge 0
\eeqq
\beqq\label{remainderterm2}
U^n_t(h) \;:=\;  \int_0^{tn} \left( \frac{S(\vth+n^{-1/2}h,s)-S(\vth,s) -  n^{-1/2} h^\top \dot S(\vth,s)}{\si(\eta_s)}\right)^2 ds   
\eeqq
\beqq\label{remainderterm3}
V^n_t(h) \;:=\;  \int_0^{tn} \frac{ [ S(\vth+n^{-1/2}h,s)-S(\vth,s) -  n^{-1/2} h^\top \dot S(\vth,s) ] [ n^{-1/2} h^\top \dot S(\vth,s) ] }{\si^2(\eta_s)}\; ds   \;. 
\eeqq
We shall show convergence of $I^n_T(\vth)$ in (\ref{leadingterm2}) to the Fisher information 
$$
I^{(\vth)} \;=\; \int_0^T \dot S(\vth,s)\,  \dot S(\vth,s)^{\!\top} \la^{(\vth)}(ds)  
$$
$Q^\vth$-almost surely  as $\nto$, weak convergence of $\call\left( \Delta^n_T(\vth) \mid Q^\vth \right)$ in (\ref{leadingterm1}) to $\caln(0,I^{(\vth)})$, as well as  
\beqq\label{aimofproof}
N^n_T(h_n) \,,\, U^n_T(h_n) \,,\, V^n_T(h_n)  \quad\lra\quad 0 \quad\mbox{in $Q^\vth$-probability as $\nto$}
\eeqq
for arbitrary bounded sequences $(h_n)_n$ in $\bbr^d$.

2)~ Consider first  the processes   $(I^n_t(\vth))_{t\ge 0}$ defined in (\ref{leadingterm2}). Write $F_{j,j'}:C_T\to\bbr$ for the function 
$$
\al  \;\lra\;  F_{j,j'}(\al) := \int_0^{T} \left( \frac{ \dot S_j(\vth,s)  \dot S_{j'}(\vth,s)}{\si^2(\al(s))} \right) ds  
$$
with arbitrary $1\le j,j' \le d$. By definition of $\la^{(\vth)}(ds)$ in (\ref{measureslambda}) and since $\dot S_j(\vth,\cdot) \in \mathbb{H}^{(\vth)}$  by $(H3)$,  
\beqq\label{leadingterm3}
m^{(\vth)}(F_{j,j'}) \;=\; 
\int_{C_T} F_{j,j'}(\al)\, m^{(\vth)}(d\al)  \;=\;  \int_0^T \dot S_j(\vth,s)  \dot S_{j'}(\vth,s)\, \la^{(\vth)}(ds)   
\eeqq
is well defined and finite, and equals   $\left[I^{(\vth)}\right]_{j,j'}$. Combining (\ref{leadingterm2}) with theorem A ,  we obtain 
\beqq\label{leadingterm4}
\left[ I^n_T(\vth) \right]_{j,\ell}  \;=\;  \frac1n \sum_{\ell=1}^{n} F_{j,j'}(\mathbb{X}_\ell) \;\lra\; m^{(\vth)}(F_{j,j'}) \;=\; \left[I^{(\vth)}\right]_{j,j'} 
\eeqq
almost surely as $\nto$. This is $Q^\vth$-almost sure convergence of $I^n_T(\vth)$ to the Fisher information $I^{(\vth)}$.

3)~ We consider the processes in (\ref{leadingterm1}):  $(\Delta^n_t(\theta))_{t\ge 0}$ are martingales with angle brackett $(I^n_t(\theta))_{t\ge 0}$ relative to $Q^\vth$ and  $(\calg_{nt})_{t\ge 0}$. 
For any $\beta\in\bbr^d$ fixed,  the function $\,F(\al) = \sum_{j,j'} \beta_j F_{j,j'}(\al)\, \beta_{j'}\,$ is nonnegative and $m^{(\vth)}$-integrable  by step 1) above, we  have 
$$
 \frac1n \sum_{k=1}^{ \lfloor \frac{tn}{T} \rfloor } F(\mathbb{X}_k) 
 \;\le\;   \beta^\top  I^n_t(\vth) \, \beta   
 \;\le\;   \frac1n \sum_{k=1}^{ \lceil \frac{tn}{T} \rceil } F(\mathbb{X}_k) 
$$
and thus by theorem A,  for every $u\ge 0$ fixed,  
\beqq\label{leadingterm5}
\lim_{\nto}\; \beta^\top  I^n_{uT}(\vth) \, \beta  \;\;=\;\;   u\cdot  \int_0^T \left( \beta^\top \dot S(\vth,s) \right)^2 \la^{(\vth)}(ds) 
\;\;=\;\; u \,\cdot\,  \beta^\top I^{(\vth)}\, \beta
\eeqq
almost surely as $\nto$. From (\ref{leadingterm5}), the martingale theorem [JS 87, VIII.3.22] gives weak convergence in the Skorohod path space $\mathbb{D}$ as $\nto$ 
$$
\call\left(  \left( \beta^\top  \Delta^n_{uT}(\vth)  \right)_{u\ge 0}  \mid Q^\vth \right)   \;\lra\;  
\call\left(   \left( \beta^\top  \wt B_u(\vth)  \right)_{u\ge 0} \right)  
$$
where $\wt B(\vth)$ is $d$-dimensional Brownian motion with covariance matrix $I^{(\vth)}$. But $\beta\in\bbr^d$ was arbitrary. So we have in particular weak convergence in $\bbr^d$ as $\nto$ 
$$
\call\left(  \Delta^n_T(\vth)  \mid Q^\vth \right)   \;\lra\;  \caln \left( 0 , I^{(\vth)} \right)  \;. 
$$
All assertions concerning convergence of the leading terms (\ref{leadingterm0}) are now proved. 

4)~  We turn to the remainder terms in (\ref{remainderterm0}) and fix any bounded sequence $(h_n)_n$ in $\bbr^d$. We start with processes  (\ref{remainderterm2}) and write, for $N$ arbitrarily large and fixed, 
$$
U^n_{NT}(h_n) \;=\; \sum_{k=1}^{nN} G(\mathbb{X}_k,n,h_n)
$$
where $G(\cdot,n,h):C_T\to [0,\infty)$ is the function 
$$
\al \;\lra\;  G(\al,n,h) :=  \int_0^{T} \left( \frac{S(\vth+n^{-1/2}h,s)-S(\vth,s) -  n^{-1/2} h^\top \dot S(\vth,s)}{\si(\al(s))}\right)^2 ds  \;. 
$$
Define with $ \rho_\vep(\vth,\cdot) \in \mathbb{H}^{(\vth)}$ introduced in $(H3)$ 
$$
G_\vep(\al) \;:=\; \int_0^T \left( \frac{ \rho_\vep(\vth,s) }{ \si(\al(s)) } \right)^2 ds \quad,\quad  \al\in C_T \,,\; \vep >0 \;. 
$$
For $(h_n)_n$ bounded by $c$, $(H3)$ allows to find for every $\vep>0$ some $n_0$ such that 
\beqq\label{hilfsbound1}
G(\al,n,h_n)   \;\;\le\;\;   c^2\, n^{-1}\, G_\vep(\al)  \quad\mbox{for all}\; \al\in C_T \;\mbox{and all}\; n\ge n_0  \;. 
\eeqq
For any $\vep>0$ fixed, the strong law of large numbers in theorem A combined with (\ref{hilfsbound1}) shows that 
$$
\limsup_{\nto}\; U^n_{NT}(h_n) \;\;\le\;\; \lim_{\nto}\;  c^2\, \frac1{n}\, \sum_{k=1}^{nN} G_\vep(\mathbb{X}_k)  \;=\;   c^2\, N\, m^{(\vth)}(G_\vep)  
$$
almost surely. But the second part of assumption $(H3)$ shows that  
$$
m^{(\vth)}(G_\vep) \;=\; 
\int_{C_T} G_\vep(\al)\, m^{(\vth)}(d\al)  \;=\;  \int_0^T \rho^2_\vep(\vth,s)\, \la^{(\vth)}(ds) \quad\lra\quad 0  
$$
as $\vep$ decreases to $0$. Both last assertions combined give 
\beqq\label{negligible1}
\lim_{\nto}\; U^n_{NT}(h_n) \;\;=\;\;  0 \quad\mbox{almost surely, for arbitrary fixed $N$ fixed} \;. 
\eeqq
The processes $(U^n_t(h_n))_{t\ge 0}$ being increasing, this implies the following property of the paths:  
\beqq\label{pathwisenegligible1}
\mbox{$U^n_\bullet(h_n)$ vanish uniformly on compact time intervals, $Q^\vth$-almost surely as $\nto$} \;. 
\eeqq

5)~ Next, $(N^n_t(h_n))_{t\ge 0}$ defined in (\ref{remainderterm1}) are martingales with angle brackett $(U^n_t(h_n))_{t\ge 0}$ relative to $Q^\vth$ and  $(\calg_{nt})_{t\ge 0}$. Hence for every $n$, the process 
$$
\left(\sup\limits_{0\le s\le t} |N^n_s(h_n)|^2\right)_{t\ge 0}
$$
is dominated in the sense of [JS 87, p.\ 35] by $\left(U^n_t(h_n)\right)_{t\ge 0}$.  As a consequence of step 2) this implies 
\beqq\label{pathwisenegligible2}
\mbox{$N^n_\bullet(h_n)$ vanish uniformly on compact time intervals, in $Q^\vth$-probability as $\nto$} \;. 
\eeqq
 
6)~ Finally, the  reasoning of step 4) combined with Schwarz inequality and step 2) shows that 
\beqq\label{pathwisenegligible3}
\mbox{$V^n_\bullet(h_n)$ vanish uniformly on compact time intervals, $Q^\vth$-almost surely as $\nto$}\;. 
\eeqq
Hence all remainder terms in (\ref{aimofproof}) vanish, and the proof of theorem 2.1 is complete.\halmos\\

{\bf 2.2 Proof of theorem 1.2 : } We assume $(H1)$, $(H3){+}(H4)$, $(H5){+}(H6)$. The proof follows the approach of Davies  [D 85, p.\ 849].

1)~ In a first step, we shall show that for every $\vth\in\Theta$ and every $c<\infty$ we have approximations 
\beqq\label{approximation1}
\sup_{\zeta\in S\cap\Theta \,,\, |\zeta-\vth|\le c\, n^{-1/2}}\; \left| \Delta_n^{(\zeta)} - \left( \Delta_n^{(\vth)} - I^{(\vth)}[n^{1/2}(\zeta-\vth)] \right) \right|  \quad=\quad   o_{Q^\vth}(1) \quad\mbox{as $\nto$} 
\eeqq 
with $S := \{ \al 2^{-k} : k\in\bbn_0\,,\, \al\in\bbz^d \}$ the set of dyadic numbers. 

From (\ref{likelihood}) and theorem 1.1 we have for points $\zeta\in S\cap\Theta$ such that $|\zeta-\vth|\le c\, n^{-1/2}$
\beam
\Delta_n^{(\zeta)} &=& \frac{1}{\sqrt{n}} \int_0^{nT} \frac{\dot S(\zeta,s)}{\si^2(\eta_s)}\, dm^{(\zeta)}_s \nonumber\\ 
&=& \frac{1}{\sqrt{n}} \int_0^{nT} \frac{\dot S(\zeta,s)}{\si^2(\eta_s)}\, dm^{(\vth)}_s 
\;-\; \frac{1}{\sqrt{n}} \int_0^{nT} \frac{\dot S(\zeta,s)}{\si^2(\eta_s)}\, \left( S(\zeta,s)-S(\vth,s) \right) ds 
\label{decomposition1}
\eeam
according to Girsanov theorem. The $Q^\vth$-martingale term on the right hand side of (\ref{decomposition1}) is 
$$
\frac{1}{\sqrt{n}} \int_0^{nT} \frac{\dot S(\zeta,s)}{\si^2(\eta_s)}\, dm^{(\vth)}_s \;\;=\;\; \Delta_n^{(\vth)} \;+\; o_{Q^\vth}(1) 
$$ 
as $\nto$: this follows if we apply $(H5)$ and theorem A to angle bracketts of 
$$
\frac{1}{\sqrt{n}} \int_0^{nT} \frac{\dot S(\zeta,s) - \dot S(\vth,s)}{\si^2(\eta_s)}\, dm^{(\vth)}_s \quad,\quad \nto \;. 
$$

Rewrite the bounded variation term on the right hand side of (\ref{decomposition1}) as 
\beao
&&\frac{1}{n} \int_0^{nT} \frac{\dot S(\vth,s)\dot S(\vth,s)^\top}{\si^2(\eta_s)}\, [\sqrt{n}(\zeta-\vth)]\, ds \\
&&+\quad \frac{1}{n} \int_0^{nT} \frac{(\dot S(\zeta,s)-\dot S(\vth,s))\dot S(\vth,s)^\top}{\si^2(\eta_s)}\, [\sqrt{n}(\zeta-\vth)]\, ds \\
&&+\quad \frac{1}{n} \int_0^{nT} \frac{\dot S(\vth,s)}{\si^2(\eta_s)}\, |\sqrt{n}(\zeta-\vth)|\, \frac{ S(\zeta,s)-S(\vth,s) - (\zeta-\vth)^\top \dot S(\vth,s) }{ |\zeta-\vth| } ds \\
&&+\quad \frac{1}{n} \int_0^{nT} \frac{\dot S(\zeta,s)-\dot S(\vth,s)}{\si^2(\eta_s)}\, |\sqrt{n}(\zeta-\vth)|\, \frac{ S(\zeta,s)-S(\vth,s) - (\zeta-\vth)^\top \dot S(\vth,s) }{ |\zeta-\vth| } ds \;. 
\eeao
The first summand equals 
$$
I^{(\vth)}\, [\sqrt{n}(\zeta-\vth)] \;+\; o_{Q^\vth}(1)
$$
by definition of the Fisher information and theorem A. The second summand, by Cauchy-Schwarz combined with $(H5)$, (\ref{fisherinformation}), theorem A and $|\zeta-\vth|\le c\, n^{-1/2}$, vanishes in $Q^\vth$-probability as $\nto$. The third summand vanishes in $Q^\vth$-probability as $\nto$ by Cauchy-Schwarz combined with $(H3)$, (\ref{fisherinformation}), theorem A and $|\zeta-\vth|\le c\, n^{-1/2}$. 
The forth summand behaves similiarly: we make use of $(H5)$, $(H3)$, theorem A and $|\zeta-\vth|\le c\, n^{-1/2}$. 
Hence we have proved the decomposition (\ref{approximation1}). 

2)~  Exactly as in [D 85], decomposition (\ref{approximation1}) combined with assumption $(H6)$ concludes the proof of   theorem 1.2. Thanks to $(H6)+(H4)$ and $\sqrt{n}$-consistency of the discretized estimator sequence $(\wt{\wt\vth}_n)_n$, we have an information with estimated parameter such that 
$$
I^{( \wt{\wt\vth}_n)} \;=\; I^{(\vth)} \;+\; o_{Q^\vth}(1) \quad\mbox{and}\quad 
(I^{( \wt{\wt\vth}_n)})^{-1} \;=\; (I^{(\vth)})^{-1} \;+\; o_{Q^\vth}(1) \;. 
$$
Exploiting (\ref{approximation1}) and $\sqrt{n}$-consistency of  $(\wt{\wt\vth}_n)_n$, we have a score with estimated parameter satisfying 
$$
\Delta_n^{(\wt{\wt\vth}_n)} \;\;=\;\; \Delta_n^{(\vth)} \;-\; I^{(\vth)}[n^{1/2}(\wt{\wt\vth}_n-\vth)] \;+\; o_{Q^\vth}(1) \;. 
$$
Hence by definition of 
$$
\vth_n^*  \;\;= \;\;  \wt{\wt \vth}_n  \;+ \; n^{-1/2}\, (I^{(\wt{\wt \vth}_n)})^{-1}\, \Delta_n^{(\wt{\wt \vth}_n)} 
$$
we have immediately from  (\ref{approximation1})
\beao
\sqrt{n}(\vth_n^*-\vth)  &=&  \sqrt{n}(\wt{\wt \vth}_n-\vth) \;+ \; (I^{(\wt{\wt \vth}_n)})^{-1}\, \Delta_n^{(\wt{\wt \vth}_n)} \\
&=&  \sqrt{n}(\wt{\wt \vth}_n-\vth) \;+ \; (I^{(\wt{\wt \vth}_n)})^{-1}\, \left( \Delta_n^{(\vth)} \;-\; I^{(\vth)}[n^{1/2}(\wt{\wt\vth}_n-\vth)] \;+\; o_{Q^\vth}(1) \right) \\
&=&  (I^{(\vth)})^{-1}\,\Delta_n^{(\vth)} \;+\; o_{Q^\vth}(1)  \quad=\quad Z_n^{(\vth)} \;+\; o_{Q^\vth}(1) 
\eeao
for the central sequence defined in (\ref{defcentralsequence}). An application of (\ref{criterionattheta}) in corollary 1.1' finishes the proof.\halmos\\

\section*{3. A minimum distance estimator sequence} 

We consider a sequence of minimum distance estimators (MDE) in the sense of Millar ([M 81], see also [Ku 04, Chapter 2.2], [H 08, Chapter 2.B+C]). We write $H$ for the Hilbert space $\,L^2( [0,T], \calb([0,T]),  \lala )\,$ with Lebesgue measure $\lala$  and assume in this section 
\beqq\label{boundedness-sigmasquared}
\mbox{$\si^2(\cdot)$ is bounded and bounded away from $0$}  
\eeqq
together with $(H1)$+$(H3)$; slighly modified analogues to $(H2)$ and $(H4)$ will come in below. 
 
To a path segment $\;\mathbb{X}_k = \left(\; \xi_{(k-1)T+s} \;\right)_{0\le s\le T}\;$ we  associate 
$$
\mathbb{Y}_k \;:=\; \left(\; \xi_{(k-1)T+s} - \xi_{(k-1)T} - \int_{(k-1)T}^{(k-1)T+s} b(\xi_r)\, dr \;\right)_{0\le s\le T} 
$$
and
$$
\mathbb{V}_k \;:=\; \left(\; \int_{(k-1)T}^{(k-1)T+s} \si(\xi_r)\, dW_r \;\right)_{0\le s\le T} 
$$
for  $k\ge 1$. Based on observation of a trajectory of the process (\ref{process}) with unknown $\vth$ up to time $nT$, we define empirical quantities  
$$
\wh \Psi_n \;:=\; \frac1n\; \sum_{k=1}^n\; \mathbb{Y}_k  \quad,\quad n\ge 1 
$$
and compare these --thanks to $T$-periodicity of the signals-- to  
$$
\Psi_\zeta \;:=\; \left(\; \int_0^s dv\, S(\zeta,v)  \;\right) _{0\le s\le T} 
$$
for all possible values of the parameter $\zeta\in\Theta$. We have the representation 
\beqq\label{MDE-4}
\wh \Psi_n  \;\;=\;\;  \Psi_\vth  \;+\; \left(\; \frac1n\; \sum_{k=1}^n\; \int_{(k-1)T}^{(k-1)T+s} \si(\xi_r)\, dW_r \;\right)_{0\le s\le T}
\eeqq
under $Q^\vth$, and prepare on some probability space a time-changed Brownian motion
\beqq\label{MDE-1}
\mathbb{W} = \left( B\circ\Phi(s) \right)_{0\le s\le T} \quad,\quad
\Phi(s) \;:=\; \int_0^s dv  \left[ \mu^{(\vth)} P^{(\vth)}_{0,v}\; \si^2 \right] \;\;,\;\; 0\le s\le T \;. 
\eeqq
With the usual conventions, our MDE sequence for the unknown parameter $\vth\in\Theta$ will be     
\beqq\label{MDE-3}
\wt \vth_n \;:=\; \mathop{ \rm arginf}_{\zeta\in\Theta}\; \| \wh \Psi_n - \Psi_\zeta \|_H  \;\;,\;\; n\ge 1 \;. 
\eeqq
We strengthen the identifiability condition $(H2)$ to 
\beqq\label{MDE-ID}
\inf_{\zeta\in\Theta, |\zeta-\vth|>\vep}\; \|  \Psi_\zeta - \Psi_\vth \|_H \;\;>\;\; 0 
\eeqq
for all $\vth\in\Theta$, and have by means of theorem A in section 2 a weak law of large numbers 
\beqq\label{MDE-WLLN}
\| \wh \Psi_n - \Psi_\vth \|_H \;\lra\;  0 \quad\mbox{in probability under $Q^\vth$ as $\nto$} 
\eeqq
together with an integrability property  
\beqq\label{MDE-DOM}
\sup_{n\ge 1}\; E_\vth\left(\, n\left(  \wh \Psi_n(s) - \Psi_\vth(s) \right)^2 \,\right) \;\le\; C\, s 
\quad,\quad 
\lim_{\nto}\; E_\vth\left(\, n\left(  \wh \Psi_n(s) - \Psi_\vth(s) \right)^2 \,\right) \;\;=\;\; \Phi(s)  
\eeqq  
for every $0\le s\le T$.  
By finite-dimensional convergence combined with (\ref{MDE-DOM}), the processes (\ref{MDE-4}) under $Q^\vth$ converge weakly in the path space $H$ to the Gaussian process $\mathbb{W}$ defined in (\ref{MDE-1}), see e.g.\ Cremers and Kadelka [CK~86].   This implies ([M 81], [Ku 04], or [H 08, 2.13]) that the MDE sequence (\ref{MDE-3}) converges in $Q^\vth$-probability to the true value $\vth$ and is tight at rate $\sqrt{n}\,$,  for all $\vth\in\Theta$. Thus we dispose of a {\em preliminary estimator sequence for $\vth$ which is tight at the required rate}. 
Already at this stage, we can apply theorem 1.2 and obtain an easy to calculate modified estimator sequence which is efficient in the sense of the  local asymptotic minimax bound given in corollary 1.1'.  
If in addition, for all $\vth\in\Theta$, the components $\dot S _j(\vth,\cdot)$ of the derivative in $(H3)$ satisfy 
\beqq\label{MDE-D}
D_j\Psi_\vth = \left( \int_0^s \dot S_j(\vth,v)\, dv \right)_{0\le s\le T} ,\; 1\le j\le d \;\;,\;\; \mbox{are linearly independent in $H$}   
\eeqq 
there is an explicit limit law for rescaled estimation errors ([M 81], [Ku 04], or [H 08, 2.23]) of the MDE sequence. Summarizing we have 

\vskip0.5cm
{\bf 3.1 Theorem: } Under $(H1)$--$(H6)$ and (\ref{boundedness-sigmasquared})+(\ref{MDE-ID})+(\ref{MDE-D}), we have for all $\vth\in\Theta$ weak convergence  
$$
{\call}\left( \sqrt{n}(\wt \vth_n-\vth) \mid Q^\vth \right) \quad\lra\quad 
\caln\left(\, 0\,,\, \Lambda_\vth^{-1}\, \Xi_\vth\, \Lambda_\vth^{-1}\, \right)  
$$
in $\bbr^d$ as $\nto$, 
where $\Lambda_\vth$ denotes the $d{\times}d$-matrix with entries $\int_0^T D_i\Psi_\vth(s) D_j\Psi_\vth(s) ds$, and 
$\Xi_\vth$ the $d{\times}d$-matrix with entries $\,\int_0^T \int_0^T D_i\Psi_\vth(s_1) \Phi_\vth(s_1\wedge s_2) D_j\Psi_\vth(s_2) ds_1 ds_2\,$ for $\Phi_\vth$ defined in (\ref{MDE-1}). One-step modification according to theorem 1.2 transforms the MDE sequence $(\wt \vth_n)_n$ into an estimator sequence $(\vth^*_n)_n$ which is asymptotically efficient at all points $\vth\in\Theta$ in the sense of corollary 1.1'. \\

Other estimator sequences  do exist which are efficient in the sense of the local asymptotic minimax bound of corollary 1.1'. \\

{\bf 3.2 Remark : } In addition to (\ref{boundedness-sigmasquared}), let $\Theta$ be bounded and assume that 
\beqq\label{sufficientcondition4}
\ul{C}|\zeta'-\zeta|^2  \quad\le\quad  \int_0^T (S(\zeta',t)-S(\zeta,t))^2\, dt \quad\le\quad \ov{C}|\zeta'-\zeta|^2  
\quad\mbox{for all $\zeta, \zeta' \in\Theta$}  
\eeqq
for suitable constants  $0 < \ul{C} \,,\, \ov{C} < \infty$  (compare (\ref{sufficientcondition4}) to (\ref{sufficientcondition2})+(\ref{sufficientcondition3})). Then maximum likelihood (MLE) 
$$
\vth^{(*,1)}_n \;=\;\mathop{\rm argmax}\limits_{\zeta\in\ov{\Theta}}\, L_{nT}^{\zeta / \zeta_0 } 
\quad,\quad n\in\bbn  
$$
(with $\zeta_0\in\Theta$ some fixed point, and $\ov{\Theta}$ the closure of $\Theta$) and Bayes estimator (BE) sequences 
$$
\vth^{(*,2)}_n  \;:=\;   \frac{ \int_{\Theta} \,\zeta\;\, L_{nT}^{\zeta / \zeta_0 }\; d\zeta }{ \int_{\Theta}  \,L_{nT}^{\zeta / \zeta_0 }\; d\zeta } 
\quad,\quad n\in\bbn  
$$
will be efficient in the sense of corollary 1.1'. 

The main ideas can be found in Ibragimov and Khasminskii [IH 81], LeCam and Yang [LY 90], Kutoyants [K 04]. They have to be adapted to the time inhomogeneous setting with $T$-periodic semigroup which is the case here. Without going into details, we mention that under (\ref{boundedness-sigmasquared})+(\ref{sufficientcondition4}) there are bounds for Hellinger type distances and affinities of form 
\beqq\label{IHbounds-1}
E_\zeta \left(\,  \sup_{0\le s\le t}\,  L_t^{\zeta' / \zeta}  \right)   \quad\le\quad   c_0 \;+\;  c_1 \left(\left\lfloor\frac{t}{T}\right\rfloor + 1 \right)  |\zeta'-\zeta|^2      
\eeqq
\beqq\label{IHbounds-2}
E_\zeta \left(\left[\, 1 \,-\,  \left( L_t^{\zeta' / \zeta}  \right)^{1/2} \right]^2 \right)    \quad\le\quad    \sum_{j=1}^3 \;c_j \left(\left\lfloor\frac{t}{T}\right\rfloor + 1 \right)^j |\zeta'-\zeta|^{2j}       
\eeqq
\beqq\label{IHbounds-3}
E_\zeta \left(\left[\, 1 \,-\,  \left( L_t^{\zeta' / \zeta}  \right)^{1/4} \right]^4 \right)    \quad\le\quad    \sum_{j=2}^5 \;c_j \left(\left\lfloor\frac{t}{T}\right\rfloor + 1 \right)^j  |\zeta'-\zeta|^{2j}         
\eeqq
\beqq\label{IHbounds-4}
E_\zeta \left( \left[  L_t^{\zeta' / \zeta}  \right]^{1/2} \right)  \quad\le\quad   \exp\left\{\, -\, k  \left\lfloor\frac{t}{T}\right\rfloor  |\zeta'-\zeta|^2  \right\}   
\eeqq
with positive constants $c_j,k$ which do not depend on $\zeta, \zeta'\in\Theta$ or on $t\ge 0$: based on ideas from the references above, a proof similiar to [HK 10, proof of lemma 5.1] goes through with   
$$
\delta_s := \frac{S(\zeta',s)-S(\zeta,s)}{\si(\eta_s)}  \quad,\quad 
d(s) := |S(\zeta',s)-S(\zeta,s)|^2  
$$ 
$$
{\ul c}\,d(s) \;\le\; \delta_s^2 \;\le\; \ov{c}\,d(s)   \quad,\quad 
\ul{C}\, \left\lfloor\frac{t}{T}\right\rfloor\, |\zeta'-\zeta|^2 \;\le\; \int_0^t d(s)ds \;\le\; \ov{C} \left(\left\lfloor\frac{t}{T}\right\rfloor + 1 \right) |\zeta'-\zeta|^2
$$
where we exploit (\ref{boundedness-sigmasquared})+(\ref{sufficientcondition4}), with suitable constants. Combining (\ref{IHbounds-2})--(\ref{IHbounds-4}) with convergence of experiments according to theorem 1.1, one follows the approach of Ibragimov and Khasminskii [IH 81] towards convergence of MLE and BE together with their moments of arbitary order. 
One can also use LeCam's 'Third lemma' to consider these estimator sequences under contiguous alternatives. 
For the sequence of local models at a reference point $\vth$ defined according to (\ref{quadraticdecomposition}) in theorem 1.1, let $\wt\cale = \{ \wt P_u : u\in\bbr^d \}$ denote the Gaussian limit experiment with central statistics $u^*$. Then $u^*$ is the MLE in the limit experiment, and also the BE ([IH 81, p.\ 180], [K 04, p.\ 134]).  
Thus MLE and BE in the Gaussian limit experiment $\wt\cale$ are equivariant estimators, the law of their error under arbitrary $u\in\bbr$ being equal to ${\call}(u^*|\wt P_0)$. For both sequences, we thus obtain  
\beqq\label{lastequation}
\lim_{\nto}\; \sup_{|u|\le C}\; \left|\;  E_{\vth+n^{-1/2}u} \left( \ell \left( n^{1/2}(\vth^{(*,i)}_{nT}-(\vth{+}n^{-1/2}u)) \right)\right) 
\;-\;  E_{\wt P_0} \left( \ell \left( u^* \right)\right)   \;\right| \quad=\quad 0
\eeqq
for arbitrary loss functions which are continuous subconvex with polynomial majorant, and for arbitrarily large constants $C$.   Combined with corollary 1.1', (\ref{lastequation}) is efficiency of both the MLE and the BE sequence, simultaneously for a large class of loss functions, in the sense of the local asymptotic minimax bound of corollary 1.1'. \\

\newpage
{\Large\bf References} 

\vskip0.2cm
[ADR 69]\quad
Az\'ema, J., Duflo, M., Revuz, D.: 
Mesures invariantes des processus de Markov r\'ecurrents. S\'eminaire de Probabilit\'es III, 
Lecture Notes in Mathematics {\bf 88}, 24--33. Springer 1969. 


\vskip0.2cm
[CLM 06]\quad
Castillo, I., L\'evy-Leduc, C., Matias, C.: 
Exact adaptive estimation of the shape of a periodic function with unknown period corrupted by white noise. 
Mathe.\ Meth.\ Statist.\ {\bf 15}, 1--30 (2006).

\vskip0.2cm 
[CK 86]\quad
Cremers, H., Kadelka, D.: 
On weak convergence of integral functions of stochastic processes with application to processes taking paths in $L^E_p$. 
Stoch.\ Proc.\ Applic.\ {\bf 21}, 305--317 (1986). 


\vskip0.2cm 
[D 85]\quad
Davies, R.: Asymptotic inference when the amount of information is random. In: Le Cam, L., Olshen, R. (Eds): 
Proc. of the Berkeley Symposium  in honour of J. Neyman and J. Kiefer. Vol. II. Wadsworth 1985. 


 

\vskip0.2cm 
[H 70]\quad 
H\'ajek, J.: A characterization theorem of limiting distributions for regular estimators.\\  
Zeitschr.\ Wahrscheinlichkeitstheor.\ Verw.\ Geb.\ {\bf 14}, 323--330, 1970. 
 

\vskip0.2cm 
[H 08]\quad 
H\"opfner, R.: Asymptotische Statistik. Manuscript in progress. Mainz 2008. \\
Under {\small\tt http://www.mathematik.uni-mainz.de/$\sim$hoepfner}

\vskip0.2cm 
[HK 10]\quad 
H\"opfner, R., Kutoyants, Y.:  
Estimating discontinuous periodic signals in a time inhomogeneous diffusion. 
Submitted to Statist.\ Inference Stoch.\ Proc., {\tt arXiv:0903.5061}. 

\vskip0.2cm
[IH 81]\quad
Ibragimov, I., Has'minskii, R.: Statistical estimation. Springer 1981. 



\vskip0.2cm
[JS 87]\quad
Jacod, J., Shiryaev, A.: Limit theorems for stochastic processes. Springer 1987. 



\vskip0.2cm
[KS 91]\quad
Karatzas, J., Shreve, S.: Brownian motion and stochastic calculus. 2nd ed.\ Springer 1991. 

\vskip0.2cm
[K 04]\quad
Kutoyants, Y.: Statistical inference for ergodic diffusion processes. Springer 2004. 


\vskip0.2cm 
[L 68]\quad
Le Cam, L.: Th\'eorie asymptotique de la d\'ecision statistique. Montr\'eal 1969. 
 
\vskip0.2cm 
[LY 90]\quad
Le Cam, L., Yang, G.: Asymptotics in statistics. Some basic concepts. 
Springer 1990.   (2nd Ed.\ Springer 2002).

\vskip0.2cm
[LS 81]\quad
Liptser, R., Shiryaev, A.: Statistics of random processes. Vols.\ I+II,  Springer 1981, 2nd Ed.\ 2001. 

\vskip0.2cm
[MT 93]\quad
Meyn, S., Tweedie, R.: Markov chains and stochastic stability. Springer 1993.

\vskip0.2cm
[M 83]\quad
Millar, P.: 
The minimax principle in in asymptotic statistical theory. 
In: P.\ Hennequin (Ed.), Ecole d'Et\'e de Probabilit\'es de St.\ Flour XI 1981.  
Lect.\ Notes in Math.\ {976}, Springer 1983. 

 

\vskip0.2cm
[R 75]\quad
Revuz, D.: Markov chains. North Holland 1975.  

\vskip0.2cm
[RY 91]\quad
Revuz, D., Yor, M.: 
Continuous martingales and Brownian motion. Springer 1991. 


\vskip0.2cm
[S 65]\quad
Skorokhod, A.: Studies in the theory of random processes. Addison-Wesley 1965.  





\vskip1.5cm
\small
Yury A.\ Kutoyants\\
Laboratoire de Statistique et Processus, Universit\'e du Maine, F--72085 Le Mans Cedex 9\\ 
{\tt kutoyants@univ-lemans.fr}\\
{\tt http://www.univ-lemans.fr/sciences/statist/pages$\_$persos/kuto.html}

\vskip0.5cm
Reinhard H\"opfner\\
Institut f\"ur Mathematik, Universit\"at Mainz, D--55099 Mainz\\ 
{\tt hoepfner@mathematik.uni-mainz.de}\\
{\tt http://www.mathematik.uni-mainz.de/$\sim$hoepfner}

\end{document}